\documentclass[12pt,thmsa]{article}
\usepackage{amsfonts}
\usepackage{amssymb}
\newtheorem{teo}{Theorem}[section]

\newtheorem{cor}[teo]{Corollary}

\newtheorem{theorem}[teo]{Theorem}

\newtheorem{obs2}[teo]{Remark}

\newtheorem{tea}{Theorem}[subsection]

\newtheorem{no2}[teo]{Note}

\newtheorem{no3}[tea]{Note}

\newcommand{\Gal}{{\rm Gal}}
\newcommand{\Frob}{{\rm Frob }}

\newcommand{\lcm}{{\rm l.c.m.}}

\newcommand{\Q}{\mathbb{Q}}

\newcommand{\PSL}{{\rm PSL}}

\newcommand{\GU}{{\rm GU}}

\newcommand{\GL}{{\rm GL}}

\newcommand{\Image}{{\rm Image}}

\newcommand{\Z}{{\mathbb{Z}}}

\newcommand{\PSU}{{\rm PSU}}

\newcommand{\cond}{{\rm cond}}

\begin{document}
\title{{\bf On the classification of geometric families of $4$-dimensional Galois
representations
}
}
\author{Luis Dieulefait
and N\'{u}ria Vila
\\
Dept. d'\`{A}lgebra i Geometria, Universitat de Barcelona;\\
Gran Via de les Corts Catalanes 585;
08007 - Barcelona; Spain.\\
e-mail: ldieulefait@ub.edu; nuriavila@ub.edu.
\thanks{Research partially supported
by MICINN grant MTM2009-07024.
 The first
named author is also partially supported by an ICREA Academia Research
Prize.
}\\
  }
\date{\empty}

\maketitle

\begin{abstract}
 We  give a classification theorem for certain four-dimensional fa\-mi\-lies of geometric $\lambda$-adic Galois representations
 attached to a pure motive. More precisely, we consider fa\-mi\-lies attached to the  cohomology of a smooth projective variety defined over  $\Q$
 with coefficients in a quadratic
  imaginary field, non-selfdual and with four different
  Hodge-Tate weights. We prove that the image is as large as possible for almost every $\lambda$ provided
   that the family is irreducible and not induced from a family of smaller dimension. If we restrict to
    semistable families an even simpler
classification
   is given. A version of the main result is given for the case where the family is attached to
    an automorphic form.\\

\end{abstract}

\vskip -20mm

\section[Introduction]{Introduction}

In this paper we study  compatible
 families of four-dimensional Galois re\-pre\-sentations.
 We  obtain  a classification theorem for families of geometric Galois representations attached to a smooth projective variety defined over  $\Q$,
 with coefficients in a quadratic
  imaginary field, non-selfdual and with four different
  Hodge-Tate weights.
In particular we have generically large images for those families of Galois representations  which are
irreducible and not induced from a family of smaller dimension.\\
This result improves our  results of \cite{DV4} where generically large image for 4-dimensional
families of Galois representations are obtained,
assuming that certain explicit conditions are satisfied.

 The first large image result of this kind was obtained by Serre for the case of
 elliptic curves without Complex Multiplication (cf. \cite{S1}). Serre also
 obtained a large image result for the four-dimensional Galois representations
 attached to principally polarized abelian surfaces with trivial endomorphism ring (cf. \cite{S2}), in his
 result he deals only with symplectic representations, while we will
 work with four-dimensional representations which are not symplectic.
  Results of generically large images
 were  obtained by the authors in \cite{D-V} for compatible families of
 three-dimensional Galois representations.\\

\bf{The Setup: Precise definition of the kind of Galois
representations that we are we going to
study.} \rm \\
\newline
Consider a compatible family of $4$-dimensional Galois
representations

$$ \rho_{\lambda} : G_{\mathbb{Q}} \rightarrow
\GL(4, K_\lambda), $$

that appear as subquotient of a family of $\ell$-adic Galois representations
  attached to the Galois action on some \'{e}tale cohomology groups $H^n_{et} (X_{\bar{\Q}}, \Q_\ell)$ of a
  smooth projective variety $X$ defined over $\Q$. We say that the family $\{ \rho_\lambda\}$ is a geometric
four-dimensional compatible family of Galois
representations of $ G_{\mathbb{Q}}$.
 We assume  that $K$ is  a quadratic imaginary field and that the Hodge-Tate numbers are $\{0,1,2,3 \}$,
  with the simplifying assumption that
  the determinants are just $\chi^6$, where  $\chi$ is the cyclotomic character.
For any prime $\lambda$ in $K$ we will
  always call $\ell$ the rational prime below it.\\
  Let $N$ denote the product of the  bad reduction primes of the variety $X$. As is well-known such a
   compatible family has a ramification set contained in the set of prime divisors of $N$ and is pure,
   i.e., the four roots of the characteristic polynomial of the image $\rho_{\lambda}(\Frob \; p)$, for
    any unramified place $p$ prime to $\lambda$,  have the same absolute value.

 From the compacity of $G_\Q$ and the continuity of the representations $\rho_{ \lambda}$,
 it follows that  we can assume that the images are contained in $\GL(4, \mathcal{O}_\lambda)$,
 where $\mathcal{O}$ denotes the ring of integers of $K$.
  This implies that we can consider the residual representations $\bar{\rho}_{ \lambda} $ with values in
 $\GL(4, \mathbb{F}_\lambda)$, obtained by composing $ \rho_{ \lambda} $ with the naive map
  ``reduction mo\-du\-lo $\lambda$".
  $$ \overline{\rho}_{\lambda} : G_{\mathbb{Q}} \rightarrow \GL(4, \mathbb{F}_\lambda).$$

Since this reduction is only well-defined up to semi-simplification,
we consider its semi-simplification, which we are going to denote
$\overline{\rho}_{\lambda}$ by abuse of notation. Thus, from now on
all residual representations will be assumed to be semi-simple.\\
\newline
Observe that from purity  it follows that if we denote by $a_p$ the trace
of the characteristic polynomial of $\rho_{ \lambda}( \Frob \; p)$ for any  prime $p\nmid \ell N$, it
holds that $|a_p| \leq 4 p^{3/2}$.\\

The characteristic polynomials of $\rho_{\lambda} (\Frob \; p)$ ($p \nmid \ell N$) will be
denoted
$$ x^4 - a_p x^3 + b_p x^2 - p^3 \bar{a}_p x + p^6 $$

The particular form of these polynomials and the fact that $b_p\in \mathbb{Z} $
follows, as is well known, from purity. We assume that the values $a_p$ generate $K$.\\

Remark: We require that the traces at Frobenius elements generate $K$, an imaginary
 quadratic field, because in the opposite case, if all traces are real, the representations
  would be selfdual, thus for all primes $\lambda$ the images would not be ``as large as possible"
   (in fact, the images would be contained in an orthogonal or a symplectic group). Observe that even with the
   assumption that the traces generate $K$, the particular form of the characteristic polynomials,
   together with Cebotarev
 density theorem imply that for any prime $\lambda$ inert in $K$ the image of $\rho_{\lambda}$ will
 be contained in the unitary group $\GU (4, \mathbb{Z}_\ell )$.\\

By conductor of the representations
$\rho_{\lambda}$, $\cond(\rho_ \lambda)$,
 we mean the prime-to-$\ell$ part of the Artin
conductor.
It is known that the conductor of the representations
$\rho_{ \lambda}$  is bounded independently of $\ell$, because the representations are geometric.
This follows from
 Deligne's corollary to  the results of de Jong (cf. \cite{Be}, Proposition
6.3.2).
We will denote by
 $c:=\lcm\{\cond(\rho_{ \lambda})\}_{\lambda}$
   this finite  uniform bound, and we will call it
    ``conductor of the family of representations".

 Let    $d$ denotes the the greatest integer such that $d^2\mid c$. \\

     {\bf Condition A)}: {\it There exists a prime $p$ congruent with $-1$ modulo $d$,  such
that none of the roots of  $Q_p(x)$, the characteristic polynomial of
$\wedge^2({\rho}_{ \lambda}) (\Frob \; p)$ , is equal to $-p^3$.}
\\

Remark: If the conductor $c$ is square-free we set condition A to be
the empty condition.\\
\newline
 This finishes the description of the geometric families of
Galois representations that we are going to consider. The following
is the main theorem of this paper, and it applies to them:

\begin{teo}
\label{teo:main} Let $\{ \rho_\lambda\}$ be a geometric
four-dimensional compatible family of Galois
representations of $ G_{\mathbb{Q}}$  with
coefficients generating an imaginary quadratic field $K$,
 with
Hodge-Tate numbers $\{0,1,2,3 \}$, with  determinant  $\chi^6$ and satisfying
 Condition A). Assume that there exists a good reduction prime $p $ such that $\Q(a_p^2) = K$ and $b_p
\neq 0$.

  Then one of the following  holds:\\
\begin{enumerate}
 \item[(i) ]
 There exists $f_1$ and $f_2$  classical cuspidal modular forms such that
$\rho_{\lambda}^{s.s.}\approx \sigma_{f_1, \lambda}\oplus (\sigma_{f_2,
\lambda}\otimes\chi)$ for every $\lambda$.
\item[(ii) ]
 There exists  $ \{ \varphi_\lambda \}$  a family of
one or  two dimensional Galois representations of the absolute
Galois group of a number field $E$ that  do not extend to $G_\Q$
such that $\rho_{\lambda}= Ind_E^{\mathbb{Q}}(\varphi_\lambda)$ for
every $\lambda$.
\item[(iii) ]
 $ \Image(\rho_{\lambda})$ is as large as possible for almost every $\lambda$.

   \end{enumerate}
\end{teo}

In other words, we have generically large image for four-dimensional
geometric families of Galois representations satisfying Condition
(A), which are non-selfdual, except if the family is reducible or
induced from a family of smaller dimension (whose members are  not
restrictions of representations of $G_\Q$). For the precise
definition of image ``as large as possible" see
 the next section and in particular the statement of Theorem \ref{teo:oldmain}.\\

 The result can be improved if we restrict to the case of geometric semistable families of representations.
  Let $N$ be (as before) the product of the primes of bad reduction
  of the variety $X$ whose cohomology is the source of the
  Galois representations $\rho_\lambda$.
  By semistable we mean that for every prime $q$ in the ramification set, i.e., $q$ dividing $N$,
   and for every $\lambda$ relatively prime to $q$, the image of the inertia group at $q$ for the
    representation $\rho_\lambda$ is pro-unipotent. In the semistable case the result can be
    proved without imposing condition (A), and only two possibilities remain, namely, we have
     the following corollary:\\

\begin{cor}
\label{teo:ss} Let $\{ \rho_\lambda\}$ be a geometric
four-dimensional  compatible family of semistable Galois
representations of $ G_{\mathbb{Q}}$  with
coefficients generating an imaginary quadratic field $K$,
 with
Hodge-Tate numbers $\{0,1,2,3 \}$ and  determinant  $\chi^6$. Assume that there exists a
good reduction prime $p $ such that $\Q(a_p^2) = K$ and $b_p
\neq 0$.

  Then one of the following  holds:\\
\begin{enumerate}
 \item[(i) ]
 There exists $f_1$ and $f_2$  classical cuspidal modular forms with trivial nebentypus such that
$\rho_{\lambda}^{s.s.}\approx \sigma_{f_1, \lambda}\oplus (\sigma_{f_2,
\lambda}\otimes\chi)$ for every $\lambda$.
\item[(ii) ]
 $ \Image(\rho_{\lambda})$ is as large as possible for almost every $\lambda$.

   \end{enumerate}
\end{cor}

According to the Langlands conjectures, certain automorphic representations of $\GL_n$
should always have geometric compatible families of Galois representations attached,
meaning that the $L$-function of the automorphic form should agree with the $L$-function
 of the family of Galois representations. It is also expected that the conductor of the
 family should agree with the level of the automorphic form. Thus, let us assume that we
  have a compatible family of  four-dimensional geometric Galois representations as in the previous
  theorem (recall that in the context of this paper, geometric implies pure)
  with the extra assumption that it comes from an algebraic automorphic form $f$
  of $\GL_4$ over $\Q$. We obtain the following corollary:

\begin{cor}
\label{teo:modular}
Let $\{ \rho_{f,\lambda} \}$ be a geometric compatible family of four-dimensional Galois
representations satisfying all the assumptions of Theorem \ref{teo:main}, and assume that
 this family is attached to an automorphic form $f$ of $\GL_4$ of $\Q$. Then the image is
  as large as possible for almost every $\lambda$,
except if one of the following holds:

- (i) $f$ is a weak endoscopic lift

- (ii) $f$ is automorphically induced from an automorphic form $g$
of $\GL_2$ of a quadratic number field or from a Hecke character of
a number field (in both cases, that are not base changed from $\Q$).

\end{cor}

Remark: Case (i), which corresponds to case (i) in Theorem
\ref{teo:main}, can not occur if we assume that the form $f$ is
cuspidal.

\section[Previous results]{Previous results and tools}
In our previous paper \cite{DV4} we studied  the images of compatible families of $4$-dimensional Galois
representations  and we proved that they are
generically large, assuming that certain explicit conditions are satisfied.
 Since the proofs of the present theorems rely heavily on the analysis and results
  of \cite{DV4}, we recall its main result (cf. Theorem 4 loc cit).

\begin{theorem}
\label{teo:oldmain} Let $ \{ \rho_\lambda \}$ be a compatible family of geometric
pure $4$-dimensional Galois representations with Hodge-Tate weights $\{0,1,2,3 \}$ and
 $N$ the product of the primes in the ramification set, with
coefficients in an imaginary quadratic field $K$ and determinant $\chi^6$. Assume that
 the following  conditions
are satisfied:\\
{ Condition 2)}: \it{There exists a prime $p$, $p \nmid N$,  such
that none of the roots of  $Q_p(x)$, the characteristic polynomial of
$\wedge^2(\rho_ \lambda) (\Frob \; p)$, is a number of the form $\eta
p^i$, where $\eta$ denotes an arbitrary root of unity and $i \in
\{1,2,3,4,5
\}$. }\\
{Condition 3)}: \it{For every quadratic character $\psi$ unramified outside $N$ there
exists a prime $p \nmid N$ with $\psi(p)=-1$ and
 $ p^3 (a_p^2 + \bar{a}_p^2) \neq a_p \bar{a}_p b_p. $}\\
{ Condition 4)}: \it{For every cubic character $\phi$ unramified outside $N$ there
 exists a prime $p \nmid N$ with $\phi(p) \neq 1$ and  $ a_p^2 b_p + p^6 \neq p^3 \bar{a}_p a_p .$}\\
 {Condition 5)}: \it{For every quadratic character $\mu$ unramified outside $N$
 there exists a prime $p \nmid N$ with $\mu(p)=-1$ and  $ a_p \neq 0.$}\\
 {Condition 6)}: \it{There exists a prime $p \nmid N$ such that: $a_p \neq \pm \bar{a}_p$.}\\
{Condition 7)}: \it{There exists a prime $p \nmid N$ such that $b_p \neq 0$ and $a_p^2$
generates $K$.}\\
Then
 the image of $\rho_\lambda$ is ``as large as possible"
for almost every prime, i.e., the image of its projectivization $P(\rho_\lambda)$ satisfies:
$$ \Image(P(\rho_\lambda)) \supseteq \PSL(4, \mathbb{Z}_\ell)$$ if $\ell$ splits in $K$, and
$$\Image(P(\rho_\lambda)) = \PSU(4, \mathbb{Z}_\ell) $$ if $\ell$ is inert in $K$.
\end{theorem}

Remarks:\\
1- Observe that condition (7) is stronger than condition (6), so in fact only five conditions
 are required, conditions (2),(3),(4),(5) and (7). By the way, conditions are numbered starting
  at (2) since at the beginning of the previous paper there was a ``condition (1)" but it was
   shown that for the compatible families we are dealing with that condition was always satisfied.\\
2- The theorem is stated in loc. cit. only for the case of coefficients in  a specific imaginary
quadratic field because the only available geometric example has coefficients in it, but the proof
 given applies to the case of coefficients in any imaginary quadratic field.\\

\rm
Another main tool that we have used at several points of our arguments is the information
 on the action of the inertia at $\ell$.
As a consequence of  results of Fontaine-Messing  and of   Fontaine-Laffaille
 (cf.  \cite{F-M} and \cite{F-L})  we know  that
 the
exponents of the fundamental characters giving the action of the residual
representation $\bar{\rho}_\lambda$ on the inertia group  $I_\ell$
must agree with the Hodge-Tate numbers of $\rho_\lambda$,  more precisely we have (see \cite{DV4}):

\begin{theorem}
\label{teo:doblestar} Consider a four-dimensional crystalline
representation $\rho_\lambda$ with Hodge-Tate numbers $\{ 0, 1, 2, 3\}$.
Then,
 if $\ell >3$, the
exponents of the fundamental characters giving the action of the (tame)
inertia subgroup $I_\ell$ of the residual representation $\bar{\rho}_\lambda$
are also $\{ 0,1,2,3\}$. Let $P_I :=
\mathbb{P}(\bar{\rho}_\lambda|_{I_\ell}^{s.s.})$
be the projectivization of the image of $I_\ell$. Then $P_I$ is a
cyclic group with order greater than or equal to $\ell - 1$.
\end{theorem}

\section[Images ]{Proof of the Main Theorem}

Let us start the proof of Theorem \ref{teo:main}. It follows from Theorem \ref{teo:oldmain} that with the running
assumptions either the image is as large as possible for almost every
prime or the residual representation falls infinitely often in one of the following two cases:\\
(a) reducible with two $2$-dimensional irreducible components\\
(b) imprimitive irreducible. \\
In fact, conditions 6 and 7 in Theorem  \ref{teo:oldmain} hold by
assumption and conditions 2 to 5 in Theorem \ref{teo:oldmain} were
conditions designed ad hoc to guarantee that cases (a) and (b) only
occur finitely many times: condition 2 took care of the reducible case
with two $2$-dimensional components (see section 5.1.2 of \cite{DV4}) and
conditions 3 to 5 took care of the imprimitive irreducible cases (see
section 5.2 of loc. cit.). Compared to Theorem \ref{teo:oldmain}, the
difference is that now we are not imposing conditions 2 to 5, and this
is why we are allowing these two cases to occur infinitely
often.\\
What remains thus is to characterize intrinsically those geometric families
 such that case (a) or (b) occurs for infinitely many primes $\ell$, so as
  to see that they correspond exactly with items (i) and (ii) in the statement of Theorem \ref{teo:main}.\\

 \subsection{Case (a)}
  Let us suppose that the residual image falls in this case for infinitely many
   primes $\lambda$. The determinants of the two-dimensional   components are
$\varepsilon\chi^i$  and
$\varepsilon^{-1}\chi^j$, respectively, where $\varepsilon$ is a character of conductor dividing $d$. \\
Let us first assume that $i\neq j$, that means $i=1$ and $j=5$ or $i=2$ and $j=4$,
consider $\varphi_\lambda:= \wedge^2\rho^{ss}_\lambda$.  The
characteristic polynomial of  $\varphi_\lambda(\Frob_p)$  modulo $\lambda$
has $\varepsilon p^i$ and $\varepsilon^{-1} p^j$ as  roots, for
$p\nmid \ell N$. If infinitely many primes $\lambda$ fall in this case, by applying
 Dirichlet's principle (and taking a suitable infinite subset of the set of primes that
  fall in this case) we can assume that the exponents $i$ and $j$ and the character
  $\varepsilon$ are independent of $\lambda$. Here we are abusing notation since we
   are thinking $\varepsilon$ as a character with values on $\mathbb{C}^*$, and we
    use the fact that there are finitely many characters with conductor dividing $d$.
We conclude that also in characteristic zero $\varepsilon p^i$ and
$\varepsilon^{-1} p^j$ are roots of the characteristic polynomial of
$\varphi_\lambda(\Frob_p)$, which is a contradiction with purity. In
fact, the roots of the characteristic polynomial of
$\varphi_\lambda(\Frob_p) = \wedge^2\rho_\lambda(\Frob_p)$ have
absolute value $p^3$, as follows from the fact that the determinant
of $\rho_\lambda (Frob_p)$ is $p^6$, $\rho_\lambda$ is
$4$-dimensional and pure, thus has roots of absolute value
$p^{3/2}$, and the roots of $\varphi_\lambda(\Frob_p)$ are products
of two roots of $\rho_\lambda (Frob_p)$.
\\
\newline
So it is enough to consider the case $i=j$, that is $i=j=3$. In this case by
 Theorem \ref{teo:doblestar} we have that the
two $2$-dimensional irreducible components of the residual mod $\lambda$
representation when restricted to the
inertia group at $\ell$ are as follows:\\
$$\pmatrix{
   \chi ^{3} & *  \cr
  0 & 1  \cr}
    \quad or \quad
\pmatrix{
   \psi_2^3 & 0 \cr
   0 & \psi_2 ^{3\ell} \cr},  $$

$$\pmatrix{
   \chi  & *  \cr
  0 & 1  \cr}\otimes \chi
    \quad  or \quad
\pmatrix{
   \psi_2 & 0 \cr
   0 & \psi_2 ^{\ell} \cr} \otimes \chi$$
   respectively, where $\psi_2$ denotes the fundamental character of level $2$.\\

The determinant of one of these components is $\varepsilon \chi^3$.
 Suppose that for infinitely many of these primes $\lambda$ the
  corresponding character $\varepsilon$ is odd, i.e., it satisfies
   $\varepsilon(c) = -1$ where $c$ denotes complex conjugation.
   Then, for every prime $p$ congruent to $-1$ modulo $d$ we have
    that $\varepsilon(p) p^3 = - p^3$ is, modulo infinitely many
    primes $\lambda$, a root of the characteristic polynomial $Q_p(x)$
     of  $\varphi_\ell(\Frob_p)$. We conclude from this that the same
     holds in characteristic $0$, i.e., that $-p^3$ is a root of the
     polynomial $Q_p(x)$  for every $p$ congruent to $-1$ modulo $d$, and this contradicts condition (A).  \\
Thus, we can assume without loss of generality that the character $\varepsilon$ is even, and
 therefore that the determinants of the two two-dimensional irreducible components are odd.\\
   By Serre's conjecture, which is now a theorem (cf. \cite{D}, \cite{K-W} and \cite{Ki}),
   this implies that for every $\lambda$ that falls in this case there are two classical
   cuspidal Hecke eigenforms $f_1$ and $f_2$ such that
   $$\rho^{ss}_\lambda\equiv \sigma_{f_1, \lambda}\oplus (\sigma_{f_2,
\lambda}\otimes\chi)\pmod \lambda$$
 Moreover applying the strong version of Serre's conjecture (cf. \cite{S3}) we can say a lot
 about the level and weight of these two modular forms: using the description of the action
 of the inertia group at $\ell$ for the two $2$-dimensional components that we gave above we
  deduce that they can be taken of weight $4$ and $2$, respectively. Using also the multiplicative
   property of conductors we conclude that $f_1\in S_4(\Gamma_1(c_1))$ and
$f_2\in S_2(\Gamma_1(c_2))$ with $c_1 c_2|c$. As a consequence, we have finitely many possibilities for the
modular forms $f_1$ and $f_2$. So, if this occurs for infinitely
many primes $\lambda$, by Dirichlet principle we can assume that $f_1$
and $f_2$ are independent of $\lambda$. Then $$\rho^{ss}_\lambda\equiv
\sigma_{f_1, \lambda}\oplus (\sigma_{f_2, \lambda}\otimes\chi)\pmod
\lambda,$$ for infinitely many $\lambda$. But this implies that the characteristic polynomials
 at unramified places
are the same, since they agree residually in infinitely many characteristics.\\
 Let us explain this step in more detail: take a density $0$ set of primes $\lambda$ as a
 subset of the infinite set of reducible primes, and then compare the characteristic polynomials
  at Frobenius elements corresponding to a density $1$ set of primes disjoint from it, thus we conclude
   that the traces agree on a density $1$ set of Frobenius elements, which by Cebotarev is enough
    to guarantee that the Galois representations are isomorphic, up to semi-simplification. We conclude that:
$$\rho^{ss}_\lambda= \sigma_{f_1, \lambda}\oplus (\sigma_{f_2,
\lambda}\otimes\chi)$$ and thus that we are in case (i) of the
theorem.

\subsection{Case (b)}

Assume that  the image of the projective residual representation is
imprimitive and irreducible. At this point we divide in two cases: the residual
 representation can be induced from a character or from a two-dimensional representation.
  Let us consider the case of a character first, i.e., let us assume that for  infinitely many
primes $\lambda$
$$\rho_\lambda\equiv Ind_E^{\mathbb{Q}}(\psi)\pmod \lambda, $$

for some  $E$  number field of degree $4$ such that its Galois closure has a Galois
 group sitting in $S_4$, and $\psi$  a character of the absolute Galois group $G_E$ that can not be extended
  to a character of $G_\Q$. The character $\psi$ takes values in $\bar{\mathbb{F}}_\ell$. \\
 A priory $E$ and $\psi$ depend on the prime $\lambda$, but by Dirichlet
 principle we can assume without loss of generality that $E$  is
 independent of $\lambda$ if we know that $E$ does not ramify at $\ell$. Let us show that this is the case, at least for
  $\ell$ sufficiently large. Moreover, if $F$ is the Galois closure of $E$, let us prove that $F$ must
  be unramified at $\ell$, for $\ell $ sufficiently large.\\
  We will use the description of the image of inertia at
  $\ell$ given in
Theorem \ref{teo:doblestar}, valid for primes $\ell >3$, $\ell \nmid
N$. In what follows we will assume that $\ell \nmid N$ and that
$\ell >5$. If we consider the image of the projectivization
$\mathbb{P} (\bar{\rho}_\lambda )$, then in this imprimitive
irreducible case this image is a group $H$ that fits in an exact
sequence:
$$ 0 \rightarrow C \rightarrow H \rightarrow S \rightarrow 0 $$
where $C$ is cyclic and $S$ is certain non-trivial subgroup of
$S_4$. In particular $| S  | \leq 24$. Note that the quotient $S$ of
$H$ is the Galois group $Gal(F/\Q)$. The subgroup $C$ is the
projectivization of a diagonal subgroup in
$\GL_4(\bar{\mathbb{F}}_\ell)$ and $H$ is contained in the
normalizer of $C$. Moreover, since elements in $H \setminus C$ act
on $C$ as permutation matrices, it follows that a matrix in $C$
having four different eigenvalues will not conmute with any element
of $H \setminus C$. From this, using the description of the
projective  image inertia at $\ell$, the fact that its order grows
with $\ell$ and that it is a cyclic group, and the fact that the
four digits of the exponent of the fundamental character are
different, we see that if $\ell$ is sufficiently large this inertia
group must be contained in $C$: if not, for large $\ell$ a large but
proper subgroup of this inertia group would be contained in $C$ and
the rest in $H \setminus C$, thus it would be non-abelian. Hence,
$F$ is unramified at
$\ell$.\\
\newline
 We conclude that the
 ramification set of $E$ is contained in the set of prime divisors of $N$,
 and the degree of $E$ is also bounded. By Hermite-Minkowski, there are only
 finitely many such fields $E$, thus we can assume (by taking a suitable infinite subset
  of the set of primes that fall in this case) that $E$ is independent of $\lambda$. Moreover, we
 have certain uniformity on the characters $\psi$, even if they are only defined as
   residual characters, because their  conductors (i.e., the prime-to-$\ell$ part of
    their Artin conductors)
 must be  bounded in terms of the conductor $c$ of the family
 $\rho_\lambda$ (this is automatic, use for example the formula giving the behavior
  of conductors under base change) and there are finitely many possibilities for
 the exponents of the fundamental characters describing the ramification at the prime
  $\ell$ because of Theorem \ref{teo:doblestar}. In this situation, we can apply a result of Serre (see \cite{S1}, $\S3$, Theorem 1, 
  for the  case of coefficients in $\Z$ and \cite{R}, Theorem (MT 2), for the generalization to the general case) to
  the restriction to the Galois closure of $E$ of the given family of Galois representations,
  and conclude that residual reducibility for infinitely many primes gives reducibility of the $\lambda$-adic representations, for every $\lambda$, thus:
 $$\rho_\lambda= Ind_E^{\mathbb{Q}}(\phi), $$
for all $\lambda$, for some character $\phi$ of the absolute Galois group of $E$.
This case is included in item (ii) of the theorem.\\

Finally, we have to consider the second possibility. Assume that for infinitely many $\lambda$, we have
$$\rho_\lambda\equiv Ind_E^{\mathbb{Q}}(\sigma_\lambda)\pmod \lambda, $$
where $\sigma_\lambda$ is a two-dimensional representation of $G_E$
that is not the restriction of a two-dimensional representation of
$G_\Q$, for $E$  a quadratic extension of $\Q$. As we have shown in
the previous case, we see in the same way that for $\ell$
sufficiently large the field $E$ will be unramified at $\ell$, thus
$E$ can be assumed (applying Dirichlet principle) independent of
$\lambda$. This induced representation is by assumption irreducible,
then
$$a_p\equiv tr(\overline{\rho_\lambda}(\Frob_p)) \equiv 0\pmod \lambda,$$
for all $p\nmid N \ell$ a prime inert in $E/\Q$. Since this holds
for infinitely many primes $\lambda$, we obtain $a_p=0$, for a set
of primes $p$ of density $1/2$ (the primes inert in $E/\Q$). By
\cite{C-S}, Theorem 3.2,  this implies that $\rho_\lambda =
Ind_E^{\mathbb{Q}}(\sigma'_\lambda)$, for some two-dimensional
family of $\lambda$-adic representations $\{ \sigma'_\lambda\}$ of
$G_E$. This concludes
 the proof, since this case is also included in item (ii) of the theorem.\\

\section{Proof of the Corollaries}

\subsection{Proof of \ref{teo:ss}}
To begin the proof of this corollary, let us explain how we can get rid of condition (A).
Recall that in the proof of Theorem \ref{teo:main} this condition was used in the case of
reducible residual image with two two-dimensional irreducible components to gain control
on the Dirichlet character $\varepsilon$ in the determinant of an irreducible component.
 Let us show that in the semistable case this Dirichlet character has to be trivial, for
sufficiently large
$\ell$.\\
The character $\varepsilon$ only ramifies at primes in $N$, and its conductor is bounded
 by $d$. From this we see that the order of $\varepsilon$ is bounded independently of $\lambda$.
  By the semistability assumption, the ramification at any prime dividing $N$ of any residual
  mod $\lambda$ representation in the family is given by an $\ell$-group. We conclude that for
   sufficiently large $\ell$ the character $\varepsilon$ has to be trivial.\\
Observe that what we have just proved has two consequences: first, that we do not need condition
(A) to guarantee that in the residually reducible case the two irreducible components will be
 odd two-dimensional representations; and second, that when proving modularity of these two
 two-dimensional components as in the proof of Theorem \ref{teo:main} since the character
 $\varepsilon$ is now trivial the two modular forms $f_1$ and $f_2$ have trivial nebentypus.\\
Finally, let us explain why in the semistable case the family of representations can
not be induced from a family of smaller dimension. As we saw in the proof of theorem
 \ref{teo:main}, this case occurs if and only if the residual representations are
 induced for almost every $\lambda$. Thus, if the family were induced, part of
 the ramification of the residual mod $\lambda$ representations at the primes
 in $N$ will come, at least for almost every $\lambda$, from the induction.
 Namely, there will be a  extension of $\Q$, either quadratic or quartic such that its Galois closure has Galois group
  a subgroup of $S_4$ such that some of the primes $q$
  dividing $N$ will ramify in this extension and the ramification of the
  residual representations at $q$ will thus have a quotient of order $2$, $3$ or $4$.
   But the semistability assumption implies that the ramification at $q$ of these
    residual representation is an $\ell$-group, thus this can not happen if $\ell >3$.

\subsection{Proof of \ref{teo:modular}}

This corollary follows easily from Theorem \ref{teo:main}. In fact, if
 the four-dimensional compatible family is attached to $f$, we already
  know that the image fails to be as large as possible for almost every
   prime only in very specific cases: If the family is as in item (i) of
    Theorem \ref{teo:main} then $f$ is a weak endoscopic lift by definition.\\
If the family is as in item (ii) of the theorem, either it is induced from a
one-dimensional or from a two-dimensional family. In the first case, clearly
$f$ is automorphically induced from a Hecke character, the character such that
 the corresponding $\ell$-adic representations appear in the formula in item (ii).
 In the second case, it remains to show that being the Galois representations induced
  from a quadratic number field $E$ the form $f$ also must be automorphically induced
  from $E$.\\
We saw during the proof of Theorem \ref{teo:main} that  in such an induced case (we
are also assuming irreducibility because reducible cases are covered by item (i))
the family satisfies $a_p = 0$ for all primes $p$ inert in $E/\Q$. Thus it satisfies:
$$\rho_{f,\lambda} = \rho_{f,\lambda} \otimes \mu $$
where $\mu$ is the quadratic character of the extension $E/\Q$. By
strong multiplicity one (cf. \cite{J-S}),
this implies: $ f = f \otimes \mu$.\\
Hence, we can apply  the characterization of automorphic forms being in the image
of solvable base change given in \cite{A-C}, Chap. 3, Theorem 4.2 (b) and Lemma 6.6 and
conclude that $f$ is automorphically induced from an automorphic form  of $\GL_2$ of the
 quadratic number field $E$. This concludes the proof of the corollary.

%

%

\end{document}